\magnification=\magstep1
\font\twelvebf=cmbx12
%
%
\font\teneusm=eusm10
\font\seveneusm=eusm7
\font\fiveeusm=eusm5
\newfam\eusmfam
\textfont\eusmfam=\teneusm
\scriptfont\eusmfam=\seveneusm
\scriptscriptfont\eusmfam=\fiveeusm

\font\tenmib=cmmib10
\font\sevenmib=cmmib7
\font\fivemib=cmmib5
\newfam\mibfam
\textfont\mibfam=\tenmib
\scriptfont\mibfam=\sevenmib
\scriptscriptfont\mibfam=\fivemib

\font\tenss=cmss10
\font\sevenss=cmss8 scaled 833
\font\fivess=cmr5
\newfam\ssfam
\textfont\ssfam=\tenss
\scriptfont\ssfam=\sevenss
\scriptscriptfont\ssfam=\fivess
\def\ss{\fam\ssfam}
\thinmuskip = 2mu
\medmuskip = 2.5mu plus 1.5mu minus 2.1mu  
\thickmuskip = 4mu plus 6mu
\def\loosegraf#1\par{{%
\baselineskip=13.4pt plus 1pt \lineskiplimit=1pt \lineskip=1.3 pt
#1\par}}

%
%
\expandafter\edef\csname amssym.def\endcsname{%
       \catcode`\noexpand\@=\the\catcode`\@\space}
\catcode`\@=11
%

\def\undefine#1{\let#1\undefined}
\def\newsymbol#1#2#3#4#5{\let\next@\relax
 \ifnum#2=\@ne\let\next@\msafam@\else
 \ifnum#2=\tw@\let\next@\msbfam@\fi\fi
 \mathchardef#1="#3\next@#4#5}
\def\mathhexbox@#1#2#3{\relax
 \ifmmode\mathpalette{}{\m@th\mathchar"#1#2#3}%
 \else\leavevmode\hbox{$\m@th\mathchar"#1#2#3$}\fi}
\def\hexnumber@#1{\ifcase#1 0\or 1\or 2\or 3\or 4\or 5\or 6\or 7\or 8\or
 9\or A\or B\or C\or D\or E\or F\fi}

\font\tenmsa=msam10
\font\sevenmsa=msam7
\font\fivemsa=msam5
\newfam\msafam
\textfont\msafam=\tenmsa
\scriptfont\msafam=\sevenmsa
\scriptscriptfont\msafam=\fivemsa
\edef\msafam@{\hexnumber@\msafam}
\mathchardef\dabar@"0\msafam@39
\def\dashrightarrow{\mathrel{\dabar@\dabar@\mathchar"0\msafam@4B}}
\def\dashleftarrow{\mathrel{\mathchar"0\msafam@4C\dabar@\dabar@}}

\def\ulcorner{\delimiter"4\msafam@70\msafam@70 }
\def\urcorner{\delimiter"5\msafam@71\msafam@71 }
\def\llcorner{\delimiter"4\msafam@78\msafam@78 }
\def\lrcorner{\delimiter"5\msafam@79\msafam@79 }
\def\yen{{\mathhexbox@\msafam@55 }}
\def\checkmark{{\mathhexbox@\msafam@58 }}
\def\circledR{{\mathhexbox@\msafam@72 }}
\def\maltese{{\mathhexbox@\msafam@7A }}

\font\tenmsb=msbm10
\font\sevenmsb=msbm7
\font\fivemsb=msbm5
\newfam\msbfam
\textfont\msbfam=\tenmsb
\scriptfont\msbfam=\sevenmsb
\scriptscriptfont\msbfam=\fivemsb
\edef\msbfam@{\hexnumber@\msbfam}
\def\Bbb#1{{\fam\msbfam\relax#1}}
\def\widehat#1{\setbox\z@\hbox{$\m@th#1$}%
 \ifdim\wd\z@>\tw@ em\mathaccent"0\msbfam@5B{#1}%
 \else\mathaccent"0362{#1}\fi}
\def\widetilde#1{\setbox\z@\hbox{$\m@th#1$}%
 \ifdim\wd\z@>\tw@ em\mathaccent"0\msbfam@5D{#1}%
 \else\mathaccent"0365{#1}\fi}
\font\teneufm=eufm10
\font\seveneufm=eufm7
\font\fiveeufm=eufm5
\newfam\eufmfam
\textfont\eufmfam=\teneufm
\scriptfont\eufmfam=\seveneufm
\scriptscriptfont\eufmfam=\fiveeufm

\csname amssym.def\endcsname

\expandafter\ifx\csname pre amssym.tex at\endcsname\relax \else \endinput\fi
\expandafter\chardef\csname pre amssym.tex at\endcsname=\the\catcode`\@
\catcode`\@=11
\newsymbol\boxdot 1200
\newsymbol\boxplus 1201
\newsymbol\boxtimes 1202
\newsymbol\square 1003
\newsymbol\blacksquare 1004
\newsymbol\centerdot 1205
\newsymbol\lozenge 1006
\newsymbol\blacklozenge 1007
\newsymbol\circlearrowright 1308
\newsymbol\circlearrowleft 1309
\undefine\rightleftharpoons
\newsymbol\rightleftharpoons 130A
\newsymbol\leftrightharpoons 130B
\newsymbol\boxminus 120C
\newsymbol\Vdash 130D
\newsymbol\Vvdash 130E
\newsymbol\vDash 130F
\newsymbol\twoheadrightarrow 1310
\newsymbol\twoheadleftarrow 1311
\newsymbol\leftleftarrows 1312
\newsymbol\rightrightarrows 1313
\newsymbol\upuparrows 1314
\newsymbol\downdownarrows 1315
\newsymbol\upharpoonright 1316
 
\newsymbol\downharpoonright 1317
\newsymbol\upharpoonleft 1318
\newsymbol\downharpoonleft 1319
\newsymbol\rightarrowtail 131A
\newsymbol\leftarrowtail 131B
\newsymbol\leftrightarrows 131C
\newsymbol\rightleftarrows 131D
\newsymbol\Lsh 131E
\newsymbol\Rsh 131F
\newsymbol\rightsquigarrow 1320
\newsymbol\leftrightsquigarrow 1321
\newsymbol\looparrowleft 1322
\newsymbol\looparrowright 1323
\newsymbol\circeq 1324
\newsymbol\succsim 1325
\newsymbol\gtrsim 1326
\newsymbol\gtrapprox 1327
\newsymbol\multimap 1328
\newsymbol\therefore 1329
\newsymbol\because 132A
\newsymbol\doteqdot 132B
 
\newsymbol\triangleq 132C
\newsymbol\precsim 132D
\newsymbol\lesssim 132E
\newsymbol\lessapprox 132F
\newsymbol\eqslantless 1330
\newsymbol\eqslantgtr 1331
\newsymbol\curlyeqprec 1332
\newsymbol\curlyeqsucc 1333
\newsymbol\preccurlyeq 1334
\newsymbol\leqq 1335
\newsymbol\leqslant 1336
\newsymbol\lessgtr 1337
\newsymbol\backprime 1038
\newsymbol\risingdotseq 133A
\newsymbol\fallingdotseq 133B
\newsymbol\succcurlyeq 133C
\newsymbol\geqq 133D
\newsymbol\geqslant 133E
\newsymbol\gtrless 133F
\newsymbol\sqsubset 1340
\newsymbol\sqsupset 1341
\newsymbol\vartriangleright 1342
\newsymbol\vartriangleleft 1343
\newsymbol\trianglerighteq 1344
\newsymbol\trianglelefteq 1345
\newsymbol\bigstar 1046
\newsymbol\between 1347
\newsymbol\blacktriangledown 1048
\newsymbol\blacktriangleright 1349
\newsymbol\blacktriangleleft 134A
\newsymbol\vartriangle 134D
\newsymbol\blacktriangle 104E
\newsymbol\triangledown 104F
\newsymbol\eqcirc 1350
\newsymbol\lesseqgtr 1351
\newsymbol\gtreqless 1352
\newsymbol\lesseqqgtr 1353
\newsymbol\gtreqqless 1354
\newsymbol\Rrightarrow 1356
\newsymbol\Lleftarrow 1357
\newsymbol\veebar 1259
\newsymbol\barwedge 125A
\newsymbol\doublebarwedge 125B
\undefine\angle
\newsymbol\angle 105C
\newsymbol\measuredangle 105D
\newsymbol\sphericalangle 105E
\newsymbol\varpropto 135F
\newsymbol\smallsmile 1360
\newsymbol\smallfrown 1361
\newsymbol\Subset 1362
\newsymbol\Supset 1363
\newsymbol\Cup 1264
 
\newsymbol\Cap 1265
 
\newsymbol\curlywedge 1266
\newsymbol\curlyvee 1267
\newsymbol\leftthreetimes 1268
\newsymbol\rightthreetimes 1269
\newsymbol\subseteqq 136A
\newsymbol\supseteqq 136B
\newsymbol\bumpeq 136C
\newsymbol\Bumpeq 136D
\newsymbol\lll 136E
 
\newsymbol\ggg 136F
 
\newsymbol\circledS 1073
\newsymbol\pitchfork 1374
\newsymbol\dotplus 1275
\newsymbol\backsim 1376
\newsymbol\backsimeq 1377
\newsymbol\complement 107B
\newsymbol\intercal 127C
\newsymbol\circledcirc 127D
\newsymbol\circledast 127E
\newsymbol\circleddash 127F
\newsymbol\lvertneqq 2300
\newsymbol\gvertneqq 2301
\newsymbol\nleq 2302
\newsymbol\ngeq 2303
\newsymbol\nless 2304
\newsymbol\ngtr 2305
\newsymbol\nprec 2306
\newsymbol\nsucc 2307
\newsymbol\lneqq 2308
\newsymbol\gneqq 2309
\newsymbol\nleqslant 230A
\newsymbol\ngeqslant 230B
\newsymbol\lneq 230C
\newsymbol\gneq 230D
\newsymbol\npreceq 230E
\newsymbol\nsucceq 230F
\newsymbol\precnsim 2310
\newsymbol\succnsim 2311
\newsymbol\lnsim 2312
\newsymbol\gnsim 2313
\newsymbol\nleqq 2314
\newsymbol\ngeqq 2315
\newsymbol\precneqq 2316
\newsymbol\succneqq 2317
\newsymbol\precnapprox 2318
\newsymbol\succnapprox 2319
\newsymbol\lnapprox 231A
\newsymbol\gnapprox 231B
\newsymbol\nsim 231C
\newsymbol\ncong 231D
\newsymbol\diagup 231E
\newsymbol\diagdown 231F
\newsymbol\varsubsetneq 2320
\newsymbol\varsupsetneq 2321
\newsymbol\nsubseteqq 2322
\newsymbol\nsupseteqq 2323
\newsymbol\subsetneqq 2324
\newsymbol\supsetneqq 2325
\newsymbol\varsubsetneqq 2326
\newsymbol\varsupsetneqq 2327
\newsymbol\subsetneq 2328
\newsymbol\supsetneq 2329
\newsymbol\nsubseteq 232A
\newsymbol\nsupseteq 232B
\newsymbol\nparallel 232C
\newsymbol\nmid 232D
\newsymbol\nshortmid 232E
\newsymbol\nshortparallel 232F
\newsymbol\nvdash 2330
\newsymbol\nVdash 2331
\newsymbol\nvDash 2332
\newsymbol\nVDash 2333
\newsymbol\ntrianglerighteq 2334
\newsymbol\ntrianglelefteq 2335
\newsymbol\ntriangleleft 2336
\newsymbol\ntriangleright 2337
\newsymbol\nleftarrow 2338
\newsymbol\nrightarrow 2339
\newsymbol\nLeftarrow 233A
\newsymbol\nRightarrow 233B
\newsymbol\nLeftrightarrow 233C
\newsymbol\nleftrightarrow 233D
\newsymbol\divideontimes 223E
\newsymbol\varnothing 203F
\newsymbol\nexists 2040
\newsymbol\Finv 2060
\newsymbol\Game 2061
\newsymbol\mho 2066
\newsymbol\eth 2067
\newsymbol\eqsim 2368
\newsymbol\beth 2069
\newsymbol\gimel 206A
\newsymbol\daleth 206B
\newsymbol\lessdot 236C
\newsymbol\gtrdot 236D
\newsymbol\ltimes 226E
\newsymbol\rtimes 226F
\newsymbol\shortmid 2370
\newsymbol\shortparallel 2371
\newsymbol\smallsetminus 2272
\newsymbol\thicksim 2373
\newsymbol\thickapprox 2374
\newsymbol\approxeq 2375
\newsymbol\succapprox 2376
\newsymbol\precapprox 2377
\newsymbol\curvearrowleft 2378
\newsymbol\curvearrowright 2379
\newsymbol\digamma 207A
\newsymbol\varkappa 207B
\newsymbol\Bbbk 207C
\newsymbol\hslash 207D
\undefine\hbar
\newsymbol\hbar 207E
\newsymbol\backepsilon 237F
\catcode`\@=\csname pre amssym.tex at\endcsname


\def\cc{{\Bbb C}}
\def\pp{{\Bbb P}}
\def\eps{\varepsilon}

\def\ss{{\Bbb S}}
\def\bb{{\Bbb B}}
\def\nn{{\Bbb N}}

\magnification=\magstep1
\centerline{\twelvebf On  convergency properties of meromorphic functions
and mappings}
\bigskip\rm
\centerline{\rm S.~Ivashkovich\footnote{}{AMS subject classification:
32 D 15, Key words: meromorphic map, Rouche theorem, spherical shell,
Fatou set.
}}

\bigskip
\def\longpoints{\leaders\hbox to 0.5em{\hss.\hss}\hfill \hskip0pt}

\bigskip
\noindent
\bf 0. Introduction.
\smallskip
\smallskip\rm
This note has two purposes. First is to discuss the possible notions of
convergency
of meromorphic functions and more generally meromorphic mappings with
values in general complex spaces. Concerning the second, it was oberved
already by Cartan and Thullen in [C-T] that the study of the domains of
existence of holomorphic functions leads to the results on the sets of
their convergency (Konvergenzbereiche). In this paper this connection
(for the case of meromorphic mappings) works in both directions.

We also give an application to the study of Fatou sets of meromorphic
self-maps of compact complex surfaces.

When one runs quicly on literature around questions of convergency of
meromorphic functions and more recently of mappings (usually with values
in $\cc\pp^n$) one founds a number of unequivalent definitions of the
notion of convergency (even in the case of functions). The reason is
that different authors study a different questions and adapt their
definitions to their problems. Let us give a simple example.

When one wishes to study the developement of, say elementary functions, into
the entire series one founds the following definition to be convenient, see
 for example [R-1] p.255.

This definition lookes nicer if one states it for the series.
 Consider a sequence $\{ \lambda_n\} $ of meromorphic functions on the
 plane domaine $\Omega\subset \cc $.

\smallskip\noindent\sl
Definition 1. \it The series $\Sigma_{n=1}^{\infty }\lambda_n$
converge compactly on $\Omega $ if for every compact $K\subset\subset \Omega $
there is a number $m=m(K)$ such that $\lambda_n$ has no poles on $K$ for
$n\ge m$
 and $\Sigma_{n\ge m}\lambda_n\mid_K$ converge uniformly on $K$.
\smallskip\rm
 For the seqence $\{ f_n:=\Sigma_{k\le n}\lambda_n\} $ of partial sums this
 notion of convergency means the following. For any compact $K\subset\subset
 \Omega $ represent each $f_n$ as a some of its
 principal part on $K$ and holomorphic function, i.e. $f_n=P_n+h_n$, where
 $P_n(z)=\Sigma_{j=1}^{l_n}{c_{n,j}\over (z-a_{n,j})^{p_{n,j}}}$. Here
 $a_{n,j}$ run over
 the poles of $f_n$ contained in $K$, and the  constants $c_j$ are uniquely
 determinated by $f_n$.
\smallskip\noindent\sl
Definition 1'. \it One says that $\{ f_n\} $ compactly converges on
$\Omega $
if for any compact $K\subset\subset \Omega $ the principal parts for $K$ of
$f_n$
stabylise for $n$ sufficiently big, and $\{ h_n\} $ converge uniformly
on $K$.
\smallskip\rm

 From the point of view of this definition the sequence $\{ {1\over z-1/n}\} $
 doesn't converge in any neighborhood of zero. Hovewer it should converge
 just by common sence. So the following definition looks also natural.
 \smallskip\noindent\sl
 Definition 2. \it A sequence $\{ f_n\} $ of meromorphic functions in
 $\Omega $
 converge  uniformly on compacts in $\Omega $ if for any point $z_0\in
 \Omega $ there is a closed
 disk $\bar\Delta (z_0,\eps )$ and a natural $m=m(z_0,\eps )$ such that
 either all $f_n, n\ge m$ are holomorphic in $\bar\Delta (z_0,\eps )$ and
 uniformly converge there, or all $1/f_n, n\ge m$ are holomorphic on
 $\bar\Delta (z_0,\eps )$ and uniformly converge there.
 \smallskip\rm
 The second definition means, in other words, that $f_n$ converges uniformly
 on compacts in $\Omega $ as a holomorphic mappings from $\Omega $ to the
 Riemann sphere $\cc\pp^1$.

 In fact meromorphic functions on plain domains are just a holomorphic
 mappings into $\cc\pp^1$, and at least from the geometrical point of view
 they are not really \it meromorphic \rm objects.

 If $\Omega \ge 2$ the meromorphic functions could have the points of
 indeterminancy. No specific value can be prescribed to a meromorphic
 function at such point, and thus from analytic point of view indeterminancies
 should be exluded from the domains of convergency, see [FS-1], (where
 this point of view comes naturally from dynamical study of holomorphic
 selfmaps of $\cc\pp^2$). Hovewer, if one looks on meromorphic function
 (or mapping) as on the analytic set - its graph, one is forced to study
 the indeterminancies as the most interesting points, carrying an essential
 information about the behavior of the converging sequence.

 Concider for example a Cremona transformation of $\cc\pp^2$: $f:[z_0:z_1:z_2]
 \to [{1\over z_0}:{1\over z_1}:{1\over z_2}]$. Then the sequence of its
 iterates $\{ f^n\} $ consists of $f$ and identity. So the maximal open
 set where the family $\{ f^n\} $ is relatively compact (i.e.Fatou set of
 $f$) should be the whole $\cc\pp^2$. Hovewer if we exclude the
 indeterminancies
 (and a fortiori their preimages) the Fatou set will be $\cc\pp^2$ minus
 three lines.

 This note consists from two parts. In Part I we give several possible
 definitions of convergency of meromorphic mappings. We discuss them
 giving {\sl pro} and {\sl contra} in each case and give some statements
  withouht proofs (not to interrupt the exposition). In Part II we give
  the main results, such as Rouche principle and apriori estimate of the
  volume, and prove the  statements from Part I.
\bigskip\noindent\bf
Table of content
\medskip \sl
\noindent
\line{0. Introduction. \longpoints pp. 1-2}

\smallskip\noindent
\line{Part I. Possible notions of convergency of meromorphic mappings.
\longpoints pp. 2-8}

1.1. Strong convergency.  1.2. Weak convergency.

1.3. $\Gamma $- convergency. 1.4. Other types of convergency.

\smallskip\noindent
\line{Part II. Proofs of the statements. \longpoints pp. 8-14}

2.1. Rouche principle. 2.2. Propagation of convergency by extension.

2.3. Propagation of extension by convergency. 2.4. Apriori estimate of

volume. 2.5. Fatou sets of meromorphic self-maps of compact complex

surfaces.

\smallskip\noindent
\line{References. \longpoints pp. 15 - 16}
 \bigskip\noindent\bf
 Part I. Possible notions of convergency of meromorphic mappings.
 \medskip\rm
 \rm Let  $\Omega $ and  $X$ be complex spaces.
 All complex spaces, which we consider in this paper are supposed to be
 reduced and normal.

 We say that a sequence fo \it holomorphic \rm mappings from a complex
 space $\Omega $ to a complex space $X$ converge uniformly on compacts in
 $\Omega $ to a holomorphic mapping $f:\Omega \to X$ if for any compact
 $K\subset\subset \Omega $ there is a compact $P\subset\subset X$ such
 that $f(K), f_n(K)\subset P$ for all $n$ and $\lim_{n\to \infty }
 \sup_ {x\in K }d(f_n(x),f(x))=0$. Here $X$ is equipped with some Hermitian
 metric. Of
 course this notion of convergency doesn't depend on the choice of this
 metric.

 \medskip\noindent\sl
 1.1. Strong convergency.
 \smallskip\rm

 Recall that a meromorphic mapping
 from $\Omega $ to $X$ is defined by a holomorphic map $f:\Omega\setminus
 F\to X$ where $F$ is analytic subset of $\Omega $ of codimension at least
 two, such that
 the closure $\bar\Gamma_f$ of its graph is an analytic subset of
 the product $\Omega\times X$. This  subset, which from now on will be noted
 as
 $\Gamma_f$ (withought bar) is called the graph of the meromomorphic map
 (again
 denoted as $f$) and clearly possesess  the following properties:

 \rm (i) $\rm{\Gamma_f}$ is an irreducible analytic subset of
 $\Omega\times X$.

 (ii) The restriction $\pi \mid _{\Gamma_f} : \Gamma_f \longrightarrow
 \Omega $ of the natural projection $\pi : \Omega \times X \longrightarrow
 \Omega $
  to $\Gamma_f$ is proper, surjective and generically one to one.

 \rm This notion of meromorphicity is due to Remmert, see [Re-2], and is based
 on the (important for us) observation that meromorphic functions on $\Omega $
  are precisely the meromorphic mappings, in the sense just defined,  into
 $X=\cc\pp^1$.

 Remark also that if $\Omega \subset \cc $ then meromorphic mappings from
 $\Omega $ are in fact holomorphic. This obvious observation is important
 in the following context. Denote by $I(f)$ the smallest analytic subset
 of $\Omega $ (now $\Omega $ is multidimensional), such that $f$ is
 holomorphic on $\Omega\setminus I(f)$. Take an irreducible complex curve
 $C\subset \Omega $, which is not contained entirely in $I(f)$. Then the
 restriction $f\mid_{C\setminus I(f)}$ extends {\it holomoprhically} onto
 the whole $C$.

 The analytic set $I(f)$ is called the indeterminancy set of $f$. If
 $z_0\in I(f)$ then the image $f[z_0]$ of $z_0$ is an analytic set,
  allwayse connected, but perhaps reducible and usually of positive
  dimension. The only thing which can be sad about the dimension of
  $f[z_0]$ is that it is not more then ${\sl dim}\Omega -1$. This follows
  from the irreducibility of $\Gamma_f$.
 When one thinks about meromorphic mappings one should keep in mind
 the following couple of examples.
 \smallskip\noindent\sl
 Example 1. (Meromorphic functions) \rm A meromorphic function $f$ on the
 complex space $\Omega $ is
 given by the open covering $\{ \Omega_j\} $ and a pairs $(h_j,g_j)\in
 {\cal O}(\Omega_j), g\not\equiv 0$, such that on $\Omega_i\cap \Omega_j$
 one  has $h_ig_j=h_jg_i$. One can allwayse suppose that $h_j$ and $g_j$ have
 relatively prime germs at all their common zeroes. Let $[z_0:z_1]$        denote
 the homogeneous  coordinates on $\cc\pp^1$. Observe now that the analytic
 set $\Gamma_f:=\{ (x,[z_0:z_1])\in \Omega_j\times \cc\pp^1:h_j(x)z_0 -
 g_j(x)z_1\} $ is correctly defined and irreducible. This will be the graph
 of the meromorphic map from $\Omega $ to $\cc\pp^1$.

 Take for example $\Omega =\cc^2$ with coordinates $x=(z_0,z_1)$ and
 $f(x)={z_1\over z_0}$. This $f$ is holomorphic on $\cc^2\setminus \{ 0\} $,
 but $f[0]=\cc\pp^1$.

 \smallskip\noindent\sl
 Example 2. (Modification) \rm
 If $\Gamma_f=\{ (z_1,z_2;[w_0:w_1])\in \cc^2\times \cc\pp^1: z_0w_1=z_1w_0
 \} $ denotes the graph of the mapping $f(z)={z_1\over z_0}$, then it is
 easy to check that $\Gamma_f$ is a manifold and that $\Gamma_f\setminus
 (\{ 0\} \times \cc\pp^1)$ is biholomorphic to $\cc^2\setminus \{ 0\} $
  under the projection onto $\cc^2$. One notes this $\Gamma_f$ usually as
  $\hat\cc^2_0$ - blown-up $\cc^2$ at origin.
 \smallskip\noindent\sl
 Example 3. (Nonextendability) \rm Meromorphic mappings, say from $\cc^2
 \setminus \{ 0\} $ with values in general complex manifold are not allwayse
 extendable to zero (Hartogs theorem in the case of holomorphic functions,
 i.e. mappings into $\cc $, and E. Levi theorem in the case of meromorphic
 functions, i.e. mappins into $\cc\pp^1$). Take, for example a Hopf surface
 $X=\cc^2\setminus \{ 0\} /(z\sim 2z)$. Then the natural projection
 $\pi :\cc^2\setminus \{ 0\} \to X$ cannot be extended meromorphically to
 zero because $\lim_{z\to 0}\pi (z)=X$.

 \smallskip

Let us intruduce our first notion of convergency. Let $\{ f_n \}$ be some
sequence of
meromorphic maps from $\Omega $ into $X$.
\smallskip\noindent\bf
Definition 1.1.1. \it  We shall say that $\{f_n \}$ strongly converges
($s$-converges) on
compacts in $\Omega $ to a meromorphic map $f :\Omega \to X
$ if for any compact $K\subset \Omega $
\smallskip
\centerline{${\cal H}- \lim_{n\rightarrow \infty}\Gamma_{f_n}\cap (K\times X)
= \Gamma_f\cap (K\times X)$}

\rm
Here by ${\cal H}-\lim $ we denote the limit in the Haussdorff metric,
supposing
that both $\Omega $ and $X$ are equipped with some Hermitian metrics. Remark
that this notion of convergency doesn't depend on a choice of metrics on
$\Omega $ and $X$.

This definition is well related with the
usual notion of convergency of holomorphic mappings. Namely, in Part II we
shall  prove
the following
\smallskip\noindent\bf
Theorem 1 \sl (Rouche principle). \it Let a sequence of meromorphic
mappings $\{ f_n\} $ between normal complex spaces $\Omega $ and $X$
strongly converge on compacts in $\Omega $ to a meromorphic map $f$. Then:
\smallskip
(a) If $f$ is holomorphic then for any relatively compact open suset $D_1
\subset D$ all

restrictions $f_n\mid_{D_1}$ are holomorphic for $n$ big enough, and
$f_n \longrightarrow f$ on compacts

in $D$ in the usual sence.

(b) If $\{ f_n\} $ are holomorphic  then $f$ is also holomorphic and
$f_n \longrightarrow f$ on compacts

in $D$ in the usual sence.
\smallskip\noindent

\smallskip\rm

We shall write $f_n \longrightarrow f$ in $D$ to denote that $\{f_n \}$
strongly converges on compacts in $D$ to a meromorphic map $f$.
 The notion of
strong convergency is well related with an extension properties of
meromorphic mappings. In one direction we have the following:
\smallskip\noindent\bf
Proposition 1.1.1. \it Let $D$ be a domain in a normal Stein space  and let
$\hat D$
be its envelope of holomorphy . Further let $\{f_n \}$ be a sequence of
meromorphic maps of $\hat D$ into a complex space $X$. Suppose that
$f_n \longrightarrow f$ in $D$ and that $f$ meromorphically extends onto
$\hat D$. Then there exists an analytic subset $A$ of $\hat D$ of codimension
at lest two, such that $f_n \longrightarrow f$ on $\hat D\setminus A$.
Moreover $A\subset I(f)$.
\smallskip\rm
The proof is given in 2.2.
\smallskip\noindent\bf
Definition 1.1.2. \it Recall, that the complex space $X$ possesses
a meromorphic (holomorphic)  extension property in dimension $n$ if for any
domain $D\subset \cc^n$ any meromorphic (holomorphic)  mapping
$f : D\longrightarrow X$ extends to the meromorphic (holomorphic)
mapping $\hat f :\hat D
\longrightarrow X$ of the envelope of holomorphy $\hat D$ of $D$ into $X$.
\smallskip

\rm Denote by

$$
H^n(r) = \left\{ (z_1,...z_n) \in C^n:\Vert z'\Vert < r, \vert z_n\vert < 1
 or \Vert z'\Vert < 1, 1 - r <\vert z_n\vert <1 \right\} \eqno(1.1.1)
$$
\medskip
\noindent
a $n$-dimensional Hartogs figure. Here $z'= (z_1,...,z_{n-1})$ and $0<r<1$.
Recall, see [Sh],[Iv-3], that  complex space $X$ possesses
a meromorphic (holomorphic) extension property in dimension $n$ iff any
meromorphic (holomorphic) mapping $f:H^n(r)\to X$ extends meromorphically
(holomorphically) onto the unit polydisk $\Delta^n$. $\Vert  \cdot  \Vert $
stands here for the polydisk norm in $\cc^n$.

If $X$ possesess a hol.ext.prop. in some dimension $n\ge 2$ then $X$
possesess this property in all dimensions, see [Sh]. We shall just say that
the space $X$ possesess the hol.ext.prop. Hovewer in [Iv-5] a 3-dimensional
compact complex manifold was constructed, which posseds a mer.ext.prop. in
dimension two but there exists a meromorphic map from punctured $3$-ball
into this manifold which doesn't extend to origin.

The Proposition 1.1.1 implies immediately the following
\smallskip\noindent\bf
Corollary 1.1.2. \it If the space $X$ possesess the hol.ext.prop. then the
maximal open subset $D\subset \Omega $ where $\{ f_n\} $ conerge is
Levi-pseudoconvex. If $X$ possesess a mer.ext.prop. in dimension $n={\sl
dim}\Omega $ then the maximal open set $D\subset \Omega $ where $\{ f_n \}$
$s$-converge is equal to a Levi-pseudoconvex set minus variety.
\smallskip\rm

Vice versa, if a meromorphic mapping from the domain $D$ in Stein space
is a strong limit  of meromorphic mappings, which are defined od the
envelope $\hat D$ of $D$, then $f$ itself extends onto $\hat D$, provided
the image space $X$ posseds a pluriclosed Hermitian metric form, see
Proposition 2.1  in Part II. This is allwayse the case when $X$ is a
compact complex surface.
\smallskip
Already in the Corollary we see that the notion of $s$-convergency
is not perfect, when one wishes to study the maximal sets where the
sequence converge. The difference becomes even more striking when one
looks for the maximal set where a given family is relatively compact.
Let ${\cal F}$ = $\{f_{\lambda}:\lambda \in \Lambda \}$ be an arbitrary
family of  meromorphic mappings of space $\Omega $ into a space  $X$.
\smallskip\noindent\bf
Definition 1.1.3. \it  One says that ${\cal F}$ is (strongly) relatively
compact on the open subset $D$ of  $\Omega $ if for any sequence $\{f_n \}
\subset {\cal F}$ there exists a subsequence $\{f_{n_k} \}$ which (strongly)
converges on compacts in $D$ to some meromorphic map
$f : D \longrightarrow X$.
The maximal open subset ${\cal N}_s$ of $\Omega $ such ${\cal F}$ is
(strongly) relatively compact   on ${\cal N}_s$ we shall call the set of
 (strong) convergency of the family ${\cal F}$.

\rm
In general nothing good can be said about  sets of strong convergency of
families
of meromorphic (and also holomorphic) mappings into a non-Stein spaces.
Namely one has the following
\smallskip\noindent\sl
Example 4. \it Let $X$ be $CP^3$ blowned up in one point. Then for every
open subset $D$ of $C^2$ one can find a sequence of holomorphic mappings of
$C^2$ into $X$  with $D$ as its set of strong convergency. \rm
\smallskip\noindent
To see this, let $(z_1,z_2,z_3)$ be coordinates of an affine part of
$\cc\pp^3$. We suppose that the blown-up point is zero in this coordinates.
For $a=(a_1,
a_2)$ and $n\in \nn $ define a mapping $f_{n,a}:\cc^2\to X$ as follows:
$f_{n,a}(z_1,z_2)=(z_1-a_1,z_2-a_2,1/n)$. If one takes $A$ to be the set of
all points in $\cc^2\setminus \bar D$ with rational cordinates, then ${\cal
F}=\{ f_{n.a}:n\in \nn ,a\in A\} $ will be the family with $N_s=D$.
\smallskip\noindent\sl
1.2. Weak convergency.
\smallskip\rm
The Proposition 1.1.1 in fact shows to us how to modify the notion of
convergency to obtain the better picture for the sets of convergency.

\smallskip\noindent\bf
Definition 1.2.1. \it  We shall say that the sequence of meromorphic mappings
from $\Omega $ to $X$ weakly converge ($w$-converge) in open subset $D$ to a
meromorphic map
$f : D\longrightarrow X$ if there exists an analytic subset $A$ of $D$
of codimention at least two such that $f_n \longrightarrow f$ in
$D\setminus A$.
\smallskip
\rm We denote this fact as $f_n\rightharpoonup f$. Changing the word strong
to weak
 in the {\sl Definition 1.1.3} we obtain the notion of the set of weak
convergency of the family of meromorphic mappings {\cal F}. We shall denote it
by ${\cal N}_w$.
\smallskip
The set of weak convergency  ${\cal N}_w$ of the family ${\cal F}$ of
meromorphic mappings from $\Omega $ into $X$ is now the maximal open suset of
$\Omega $  such that ${\cal F}$ is weakly relatively compact on ${\cal N}_w$.

The following  corollary from {\sl Proposition 1.1.1}, shows the advantage
of the second definition. Recall
that open subset $D$ of complex space $\Omega $ is called Levi-pseudoconvex
if for any poind $p\in \partial D$ one can found an open neiborhood $U$
of $p$ in $\Omega $ such that $U\cap D$ is Stein. We say that a complex
space possesses a meromorphic extension property if for any domain $D$ in
Stein space any meromorphic mapping $f:D\longrightarrow X$ extends to a
meromorphic map $\hat f :\hat D\longrightarrow X$ of the envelope of
holomorphy $\hat D$ of the domain $D$ into $X$. For example any compact
K\"ahler manifold possesses a meromorphic extension property, see [Iv-3].
\smallskip\noindent\bf
Corollary 1.2.1. \it (a) Set of weak convergency alwayse  exist.
\smallskip\noindent
(b) If space $X$ possesses a mer.ext.prop. in dimension $n={\sl dim}\Omega $,
 than the set of weak
convergency is Levi-pseudoconvex for any family of meromorphic maps of
$\Omega $ into $X$.
\smallskip\noindent\sl

\rm In fact the situation is essentially better then it is described in this
statement. Namely, from our characterisation of obstructions for the
extendibility of meromrophic mappings from domains in $\cc^n$
into the spaces with pluriclosed metrics, we shall derive that if a
meromorphic map $f:D\to X$ is a limit on $D$ (weak or strong, doesn't matter)
of meromorphic maps $f_n:\hat D\to X$, then $f$ has extra extension properties.
In fact $f$ extends onto  the envelope $\hat D$, and thus $f_n\rightharpoonup
f$
on $\hat D$! See 2.3.

\smallskip\noindent\bf
Remarks 1. \rm
 Let us give one more reason, why the second definition of convergency for
meromorphic mappings is natural. Take $\cc\pp^N$ as $X$ now. Then one can
show that any meromrophic map $f:\Omega \to \cc\pp^N$ can be locally presented
as $f=[f^0:...:f^N]$ in homogeneous coordinates, where $f^j$ are holomorphic
functions. One can show that the following is true
\smallskip\noindent\bf
Proposition 1.2.1. \it $f_n\rightharpoonup f$ iff there can be find
such local presentations $f_n=[f_n^1:...:f_n^N]$ that $f_n^j\to f^j$
as holomorphic functions for all $j=1,...,N$.
\smallskip
\rm Such type of convergency
of meromorphic mapping with values in projective manifolds was considered by
Fujimoto, see [Fj].
\smallskip\noindent\bf
2. \rm Observe that for the family ${\cal F}$ constructed in {\sl Example 4}
${\cal N}_w=\cc^2$.
\smallskip\noindent
3. \rm Hovewer, we should remark that notion of weak convergency is not as
well related with convergency of holomorphic mappings as strong does. In
particilar the Rouche principle is not longer valid in this case. To see
this consider in Example 4 the sequence $f_n(z_1,z_2)=(z_1,z_2,{1\over n})$.
$\{ f_n\} $ are holomorphic, and they wakly converge to mapping $f$ from
$\cc^2$ to $X$, which has indeterminancy at zero.

\smallskip\noindent
\sl  1.3. $\Gamma $-convergency.\rm
\medskip\rm
Fix some Hermitian metric forms  $w_X$ and  $w_{\Omega }$ on $X$ and $\Omega
$ respectively. By $p_1:\Omega
\times X \longrightarrow \Omega $ and $p_2:\Omega\times X\longrightarrow
X$ we denote the projections onto the first and second factors. On the product
$\Omega\times X$ we consider the metric form $w = p^{\ast }_1w_{\Omega } +
p^{\ast }_X $.
It will be convenient sometimes for us to consider instead of mappings
$f :\Omega \longrightarrow X$ their graphs $\Gamma_f$. By $\hat f = (z,f
(z))$ we shall denote the mapping into the graph $\Gamma _f \subset\Omega
\times X$. The volume of the graph $\Gamma _f$ of the mapping
$f$ is given by
\smallskip
$$
{\sl vol}(\Gamma _f) = \int_{\Gamma _f}w^q =
\int_{\Omega }(f^{\ast }w_X + w_{\Omega })^q \eqno(1.3.1)
$$
\smallskip
\noindent
Here by $f^{\ast }w_X$ we denote the preimage of $w_X$ under $f $, i.e.
$\phi ^{\ast }w_X = (p_1)_{\ast }p^{\ast }_2w_X$.

Recall that the Hausdorff distance between two subsets $A$ and $B$ of the
metric space $(Y,\rho )$ is a number $\rho (A,B) = \inf \{ \varepsilon : A^
{\varepsilon } \supset B, B^{\varepsilon }\supset A\} $. Here by $A^
{\varepsilon }$ we denote the $\varepsilon $-neighborhood of the set $A$, i.e.
$A^{\varepsilon } = \{y\in Y: \rho (y,A) < \varepsilon \}$.

Note that if the family $\{ {\cal F}\subset Mer(\Omega ,X)\} $ is
$s$-normal on $\Omega $ then
\smallskip
\it currents $\{ (w_{\Omega }+f_{\lambda }^*w_X)^j:f_{\lambda }\in {\cal F} \}
 $ have
uniformly bounded masses on compacts

in $\Omega $ for $j=1,..,q=dim \Omega $.
\smallskip\rm
This means, in other words that the volumes of the graphs $\Gamma_{f_
{\lambda }}$ are uniformly bounded over the compacts in $\Omega $. So,
one naturally has one more notion of convergency of the sequence $\{ f_n\} $
of meromorphic mappings of the
complex space $\Omega $ into the complex space $X$.
\smallskip
\noindent
\bf Definition 1.3.1. \rm   We shall say that $\{ f_n\}$ $\Gamma $-
converge on the compacts
 in $\Omega $  if for every
relatively compact open $D_ \subset \subset \Omega $ the sequence graphs
$\Gamma _{f_n}
\cap (D\times X)$ converge in the Hausdorff metric on $D\times X$.

One has the following
\smallskip
\noindent
\bf Proposition 1.3.1. \it Let $\{f_n\}$ be a sequence of meromorphic
mappings into a
complex space $X$. Suppose that there exists a compact $K\subset X$ and a
constant $C<\infty $ such that:

a) $f_n(\Omega )\subset K$ for all $n$;

b) ${\sl vol}(\Gamma _{f_n})\le C$ for all $n$.
\smallskip
\noindent
Then there exists a subsequence $\{f_{n_j}\}$ and a proper analytic set $A
\subset \Omega $ such that:

1) the sequence $\{\Gamma _{f_{n_j}}\}$ converges in the Hausdorff metric
to the analytic subset $\Gamma $ of $\Omega\times X$ of pure dimension $q$;

2) $\Gamma = \Gamma _{\phi }\cup \hat \Gamma $,where $\Gamma _{\phi }$ is the
graph of some meromorphic mapping $f :\Omega \longrightarrow X$,and
$\hat \Gamma $ is a pure $q$-dimensional analytic subset of $\Omega\times X$
,mapped by the projection $p_1$ onto $A$;

3) $f_{n_j}\longrightarrow f $ on compacts in $\Omega\setminus A$;

4) one has
$$
\lim_{j \longrightarrow \infty }{\sl vol}(\Gamma _{f_{n_j}})\ge
{\sl vol}(\Gamma _{f }) + {\sl vol}(\hat \Gamma ).\eqno(1.3.2)
$$

5) For every $1\le p\le dimX - 1$ there exists a positive constant $\nu _p
= \nu _p(K,h)$ such that the volume of every pure $p$-dimensional compact
analytic subset of $X$ which is contained in $K$ is not less then $\nu _p$.

6) Put $\hat \Gamma = \bigcup_{p=0}^{q-1}\Gamma _p $, where $\Gamma _p$ is a
union of all irreducible components of $\hat \Gamma $ such that ${\sl dim}
[p_1(\Gamma _p)] = p$. Then
\smallskip
$$
{\sl vol}_{2q}(\hat \Gamma ) \ge \sum_{p=0}^{q-1}{\sl vol}_{2p}(A_p)\cdot
\nu _{q-p}\eqno(1.3.3)
$$
\smallskip
\noindent
where $A_p = p_1(\Gamma _p)$.
\smallskip\rm
 The proof uses the Harvey-Shiffman generalisation of Bishop's convergency
 theorem for analytic sets and can be found in [Iv-5].

 The notion of $\Gamma $-convergency seems to be \it to week \rm to reflect
 the fact the it is the mappings are converging and not just an analytic sets.
 Hovewer it might be interesting from the measure theoretic point of vew.
  Indeed the set ${\cal N}_{\Gamma }$ of $\Gamma $-convergency of the
  family ${\cal F}$ is exactly
  the maximal open set where the currents $\{ (w_{\Omega }+f^{\ast }w_X)^n:
  f\in {\cal F}\}$ have uniformly bounded masses. We shall discuss this in
  2.4.
  \medskip\noindent\sl
  1.4. Other types of convergency. \rm
  \smallskip
  One can give also other definitions of convergency of sequences of
  meromorphic functions and mappings. One is in the spirit of Definition 1
  from the Introduction. Let a sequence of meromorphic maps $f_n:\Omega
  \to X$ is given.
  \smallskip\noindent\bf
  Definition 1.4.1. \it One says that $\{ f_n\} $ converges on compacts in
  $\Omega $ if for any relatively compact open $D\subset \Omega $
  there are

  (1) a proper modification $\pi_D:\hat D\to D$;

  (2) $N$ depending on $D$;

  \noindent
  such that the pullbacks $f_n\circ \pi_D:\hat D\to X$ are holomorphic on
  $\hat D$ for $n\ge N$ and converge there as the sequence of holomorphic
  maps.
  \smallskip\rm
  Of course this is a very restrictive notion, which means that
  indeterminancies of the sequence stabylise. So, it is a direct analog of
  the first definition from the Introduction, and is reasonable in the
  context of studying of developping of entire (for example) functions into
  the series of their principal parts.

  Such convergency implies the strong one, but not vice versa. Really,
  consider a sequence of meromorphic functions $f_n(z_1,z_2)={z_1- {1\over n}
  \over z_2}$. This sequence converge in the strong sence but indeterminancies
  do not stabylise.
  \smallskip
  In the dynamical study of holomorphic and meromorphic selfmaps the
  sequence $\{ f^n\} $ of iterates of mapping $f:X\to X$ is said to be
  normal on $D\subset X$ if it is equicontinuous there, [FS-1].
  In particular this
  exludes from $D$ the indeterminancies of $f$ and their preimages. See
  example of Cremona transformation from the Introduction in this regard.

\bigskip
\noindent\bf
Part II. Proofs of the statements.
\medskip\noindent\sl
2.1. Rouche principle.

\smallskip\noindent
Proof.

\medskip\rm
Let $\{ f_n:\Omega \to X\} $ our sequence, which strongly converge
on compacts in $\Omega $ to a meromorphic map $f:\Omega \to X$.
\smallskip\noindent
(a) Suppose that $f$ is holomorphic. Take a point $a\in \Omega $ and let
$V\ni a$ and $W\ni f(a)$ are Stein neighborhoods such that ${\overline f(V)}
\subset W$. Then $\Gamma_{f\mid_{\bar V}}\cap \partial (V\times W)\subset
\partial V
\times  W$. So the natural projection $p_1\mid_{\Gamma_f\cap (V\times W)}
:\Gamma_f\cap (V\times W)\to V$ is proper and in fact bijective.

From the strong convergency of our sequence we have for $n>>0$
$\Gamma_{f_n\mid_{\bar V}}\cap \partial (V\times W)\subset \partial V
\times  W$ and thus  $p_1\mid_{\Gamma_{f_n}\cap (V\times W)}
:\Gamma_{f_n}\cap (V\times W)\to V$ is proper and surjective. Now $V\times W$
 is Stein, so doesn't contain a compact subvarieties of positive dimension.
 So $\Gamma_{f_n}\cap (V\times W)$ is a graph of (may be multivalued)
 holomorphic correspondance from $V$ to $W$. But over a dense set $V\setminus
 I(f_n)$ this corresppondance is one-valued. Thus $\Gamma_{f_n}\cap (V\times W
 )$ is a graph of a holomorphic mapping.

\smallskip\noindent
(b) Suppose now that $f$ is not holomorphic. We must prove that $f_n$ are
also not holomorphic for $n$ big enough. Take a point $a$ on $I(f)$ such
that $k:={\sl dim}f[a]\ge 1$.
Consider two cases.
\smallskip\noindent\sl
Case 1. $k<{\sl dim}X$. \rm

Put $p={\sl dim}\Omega $. Fix two distinct points $b_1,b_2\in f[a]$. Take
a nonintesecting neighborhoods $V_i\ni b_1$ and $V_2\ni b_2$ of thouse
points in $X$ and proper holomorphic embeddings $\phi_j:V_j\to \Delta^k
\times \Delta^{m_j}$ which extend to the neighborhood of $\bar V_j$ and such
that $\phi_j(V_j\cap f[a])\cap (\bar\Delta^k\times \partial \Delta^{m_j})=
\emptyset $. Take a neighborhood $W\ni a$ in $\Omega $ and consider an
embeddings $\Phi_j:W\times V_j\to W\times (\Delta^k\times \Delta^{m_j})$
given by $\Phi_j(w,v)=(w,\phi_j(v))$. Consider a projections $\pi_j:
W\times (\Delta^k\times \Delta^{m_j})\to W\times \Delta^k$. If $W$ was
choosen sufficiently small then still $\Phi_j(\Gamma_f\cap (W\times V_j))
\cap (W\times \Delta )\partial \Delta^{m_j}=\emptyset $. Thus the
restrictions $\pi\mid_{\Phi_j(\Gamma_f\cap (W\times V_j))}:
\Phi_j(\Gamma_f\cap (W\times V_j))\to W\times \Delta^k$ will remain proper.
Denote by $A_j=\pi_j(\Phi_j(\Gamma_f\cap (W\times V_j)))$ the $p$-dimensional
 analytic subsets of $W\times \Delta^k$. Remark that
 $A_1\cap A_2\ni a$ and that $dimA_1+dimA_2=2p>p+k=dim(W\times \Delta^k)$,
 while $k=dimf[a]\le p-1$. Choosing the coordinate morphisms $\phi_j $
 genericalmly, we can suppose that $A_1\not= A_2$.

 For $n$ big enouhg we have also $\Phi_j(\Gamma_{f_n}\cap (W\times V_j))
 \cap (W\times \Delta^k)\times \partial \Delta^{m_i}=\emptyset $. Thus
 $A_j^n:=\pi_j(\Phi_j(\Gamma_{f_n}\cap (W\times V_j)))$ will be an
 analytic subsets of dimension $p$ in $W\times \Delta^k$. We have by the
 definition of strong convergency, that $\Gamma_{f_n}\cap (W\times V_j))
 \to \Gamma_f\cap (W\times V_j)$ in Haussdorff metric. This impyes that
 $A_j^n\to A_j$ in Haussdorff metric in $W\times \Delta^k$.

 The usual Rouche theorem for holomorphic functions imply now that
 $A_1^n\cap A_2^n\not= \emptyset $ for $n$ big enough. This gives us a
 point $a_n\in W$ for which $f_n[a_n]$ consists of more that one point.
 I.e. $a_n$ is a point of indeterminancy of $f_n$.
 \smallskip\noindent\sl
 Case 2. $k=dimX$, i.e. $f[a]=X$.
 \smallskip\rm
  In this case take two germs of hypersurfaces $F_j\ni b_j$ in $X$. Put
  $A_j:=\pi ((F_j\times W)\cap \Gamma_f)$, where $\pi :\Omega \times X\to
  \Omega $ is a natural projection. For $b_1,b_2$ general enough $A_1$
  and $A_2$ will be different. Now repeat the reasonings of Case 1 to
  get the same conclusion.
  \smallskip
  \hfill{q.e.d.}

\medskip\noindent\sl
2.2. Propagation of convergency by extension.
\smallskip\rm
In this paragragh we shall prove
the Proposition 1.1.1 and Corollary 1.2.1 from {\sl Part I}, which
establish the main properties of the sets of
convergency of the sequences of meromorphic mappings.
\smallskip\noindent\sl
Proof of Proposition 1.1.1. \rm Let $D'$ be a maximal open subset of $\hat D$
where $f$ is holomorphic and $\{ f_n\} $ converge to $f$ on compacts as a
sequence of holomorphic maps. Suppose that there exists a point $p\in
\partial D'\setminus I(f)$ such that $D'$  is not pseudoconvex at $p$. Take
a Stein neighborhood $V\ni p$ such that the envelope of holomorphy of
 $(D'\cup V)$ contains $p$ and $V\cap I(f)=\emptyset $.

 In the product $\hat D\times X$ take a Stein neighborhood $W$ of $\Gamma_{
 f\mid_{\hat V}}$. Now we have a sequence of \it holomorhpic \rm maps
 $f_n:V\to W$ converging to $f$, which holomorphically extends onto $\hat V$.
 $W$ is Stein, so from the maximum  principle it follows that $f_n$  converge
 to $f$ on $\hat V$. This contradics the maximality of $D'$.

 This implies that $D'\supset \hat D\setminus I(f)$.
 \smallskip
 \hfill{q.e.d.}

\smallskip\noindent\sl
Proof of Corollary 1.2.1. \rm (a) We need to prove that if for open $N_1,N_2\subset \Omega $ the
 family ${\cal F}$ is weakly relatively compact on $N_i,i=1,2$, then ${\cal F}$
 is weakly relatively compact on $N_1\cup N_2$. For this take some sequence
$\{ f_n\}
\subset {\cal F}$. We find an analytic subset $A_1\subset N_1$ of codimension
at least two and substruct from from $\{ f_n\} $ a subsequence (still denoted
as $\{ f_n\} $), which converge on compacts in $N_1\setminus A_1$ to some
meromorphic map $f_1:N_1\to X$. From this
subsequence we can substruct a subsequence once more, which will converge on
compacts in $N_2\setminus A_2$ for some analytic of codimension at least two
set $A_2$ in $N_2$ to a meromorphic map $f_2:N_2\to X$.

$f_1$ and $f_2$ will clearly coinside on $(N_1\cap N_2)$ and our subsequence
will strongly converge to $f:=f_1=f_2$ on $(N_1\cup N_2)\setminus (A_1\cup
A_2)$. From {\sl Proposition 1.1.1} we have that $I(f)\subset (A_1\cup A_2)$.
If there would exists some $a\in (A_1\cup A_2)\setminus I(f)$, then again by
{\sl Proposition 1.1.1} our subsequence would converge to $f$ in the
neighborhood
of $a$. Thus $ f_{n_k}\to f$ on compacts in $(N_1\cup N_2)\setminus I(f)$.
\smallskip\noindent
(b) This is again immediate corollary of {\sl Proposition 1.1.1}.
\smallskip
\hfill{q.e.d.}

\medskip\noindent\sl
2.3. Propagation of extension by convergency.
\smallskip\rm

\rm In fact the picture from the previous paragragh can be in some cases
reversed. Namely, from our characterisation of obstructions for the
extendibility of meromorphic mappings from domains in Stein manifolds
into the spaces with pluriclosed metrics, we shall here derive that if a
meromorphic map $f:D\to X$ is a limit on $D$ (weak or strong, doesn't matter)
of meromorphic maps $f_n:\hat D\to X$, then $f$ has extra extension properties.
In fact $f$ extends onto  the envelope $\hat D$, and thus $f_n\rightharpoonup
f$
on $\hat D$!

We call a Hermitian metric form $w$ on complex space $X$ pluriclosed if
 $dd^cw=0$. One can show that the following duality holds:

\it either compact complex manifold $X$ admits a pluriclosed metric form,
 or $X$ carries

a bidimension $(2,2)$ current $T$ such that $dd^cT$ is also positive.

\rm
We shall restrict ourselves with the following class of spaces:
\smallskip\noindent\bf
Definition 2.3.1.  \it Call the complex space $X$ disk-convex if for any
compact $K\subset X$ there is a compact $\hat K\subset X$ such that for
every analytic disk $\phi :\bar \Delta \to X$ with $\phi (\partial \Delta )
\subset K$ one has $\phi (\Delta )\subset \hat K$.
\smallskip\rm
In [Iv-3] the following theorem is proved:
\smallskip\noindent\bf
Theorem 2.3.1. \it Let $f: D \longrightarrow X$ be a meromorphic mapping
from a domain in Stein manifold $\Omega $ into a disk-convex complex space
$X$ which possesses a pluriclosed Hermitian metric form  $w$. Then
$f$ extends onto $\hat D\setminus A$ where $A$ is an analytic subset
of $D$ of pure codimension two (may be empty). Moreove, if $A$ is
nonempty then for any $\ss^3\subset \hat D\setminus A$, such that $\ss^3$
 is not homologous to sero in $\hat D\setminus A$ the image $f(\ss^3)$ is
 also not homologous to zero in $X$.
\smallskip\rm

 Recall, see [Iv-2], that a spherical shell of dimension two in complex space
$X$ is an image $\Sigma $ of the standard sphere ${\bf S^3}\subset \cc^2$
under the holomorphic map of some neighborhood of ${\bf S^3}$ into $X$, such
that $\Sigma $ is not homologous to zero in $X$. This notion is close to
the notion of the global spherical shell, introduced by Kato, see [Ka].
Thus we obtain the following
\smallskip\noindent\bf
Corollary 2.3.2. \it Let $X$ be a disk-convex complex space which possess a
pluriclosed Hermitian metric form. Then the following is equivalent:

(a) $X$ possesses a meromorphic extension property in all dimensions.

(b) $X$ contains no spherical shells.
\smallskip\rm
We turn  now to the sets of convergency of families of meromorphic
mappings.

\smallskip\noindent\bf
Corollary 2.3.3. \it If $X$
admits a pluriclosed Hermitian metric, then for any family of meromorphic
mappings from $\Omega $ to $X$ the set of weak convergency is
Levi-pseudoconvex.
\smallskip\noindent\sl
Proof. \rm Let $D=N_w$ be the set of weak convergency of our family ${\cal F}$.
 Suppose that $p\in \partial D$ is not pseudoconvex boundary point. Let
 $U\ni p$ be a neighborhood such that for some component of $U\cap D$, say
 $V$,
the projection of the envelope of holomorphy of $V$ into $U$ contains a
neighborhood $W$ of $p$.

Take some sequence $\{ f_n\} \subset {\cal F}$. Then some subsequence
$\{ f_{_k}\} $ weakly converge on $D$ to a meromorphic map $f:D\to X$. This
means by definition that there is an analytic set $A$ of codimension at
least two in $D$ such that $f_{n_k}\to f$ on compacts in $D\setminus A$. By
{\sl Theorem 2.3.1} $f$ extends meromorphically to $W$ minus analytic set of
codimension two. By {\sl Proposition 1.1.1} $f_{n_k}\to f$ on $W\setminus
(A\cup B)$.

All that left to prove is that $B\subset A$, i.e. $f$ is meromorphic on the
whole $W$. Would $B\not\subset A$ then there by {\sl Theorem 2.3.1} would be a
three-sphere $\ss^3\subset W\setminus B$ such that $f(\ss^3)\not\sim 0$ in
$X$. But $f_{n_k}\to f$ on $\ss^3$ and $f_{n_k}(\ss^3)\sim 0$ in $X$ because
$f_{n_k}$ are meromorphic on the whole $W$. This is a contradiction.
\smallskip
\hfill{q.e.d.}
\medskip\noindent\sl
2.4. Apriori estimate of the volume.
\smallskip\rm

Our aim in this paragraph is the following
\smallskip\noindent\bf
Theorem 2. \it Let $\{ f_{\lambda }\} $ be a family of meromorphic mappings
from $\Delta^2$ to a disk-convex K\"ahler space $(X,w)$ such that
$f_{\lambda }$ are
holomorphic on $A^2(1/2,1):=\Delta^2(1)\setminus \bar\Delta^2(1/2)$ and
equicontinuous there. Then
$\{ {\sl vol}(\Gamma_{f_{\lambda }})\} $ are uniformly bounded.
\smallskip\noindent\sl
Proof. \rm Consider a family of currents $T_{\lambda }:=f^*w$ on $\Delta^2$.
Write $T_{\lambda }=i/2t^{\alpha \bar\beta }_{\lambda }dz_{\alpha }\wedge
d\bar z_{\beta }$, where $t^{\alpha \bar\beta }_{\lambda }$ are distributions
on $\Delta^2$, smooth on $\Delta^2\setminus S_{\lambda }$. Here $S_{\lambda }
$ is a finite set for each $\lambda $.

Consider functions
$$
\mu_{\lambda }(z_1) = i/2\cdot \int_{\Delta_{z_1}}t^{2\bar 2}_{\lambda }
dz_2d\bar z_2. \eqno(2.4.1)
$$
\noindent
$\{ \mu_{\lambda }\} $ are uniformly bounded on $A(1/2,1)$. Moreover
$$
{\partial \over \partial z_1}\mu_{\lambda } = i/2\cdot \int_{\Delta_{z_1}}
{\partial \over \partial z_1}t^{2\bar 2}_{\lambda }dz_2d\bar z_2 =
i/2\cdot \int_{\Delta_{z_1}}{\partial \over \partial z_2}t^{1\bar 2}_
{\lambda }dz_2d\bar z_2 =
$$
$$
 = i/2\cdot \int_{\partial \Delta_{z_1}}t^{1\bar 2}_{\lambda }d\bar z_2  .
\eqno(2.4.2)
$$
This gives us the boundedness of the differentials of $\{ \mu_{\lambda }\} $
on $\Delta $ and thus the boundedness in $L^1(\Delta^2)$ of $\{ t^{2\bar 2}
_{\lambda }\} $.

In the same way we get a boundedness in $L^1(\Delta^2)$ of $\{ t^{1\bar 1}
_{\lambda }\} $. From positivity of $T_{\lambda }$ we get
$$
\int_{\Delta^2}\vert t^{1\bar 2}_{\lambda }\vert \le
\int_{\Delta^2}\sqrt{t^{1\bar 1}_{\lambda }}\cdot \sqrt{t^{2\bar 2}_{\lambda }}
\le
$$
$$
 \le \sqrt{\int_{\Delta^2}t^{1\bar 1}_{\lambda }}\cdot
\sqrt{\int_{\Delta^2}t^{2\bar 2}_{\lambda }}.\eqno(2.4.3)
$$
\noindent
This gives us the boundedness of $\{ T_{\lambda }\} $ in $L^1(\Delta^2)
$.

To estimate ${\sl vol}\Gamma_{f_{\lambda }}=\int_{\Delta^2}(f^*w+dd^c\Vert z
\Vert^2)^2$ we need to estimate also $\int_{\Delta^2}(f^*w)^2$.

Take a smooth function $1\ge \eta \ge 0$ in $\Delta^2$, $\eta\mid_{\Delta^2
(3/4)}\equiv 1$, $\eta\mid_{A^2(7/8,1)}\equiv 0$ and consider the following
potentials:
$$
U_{\lambda }(z)=-\int_{\cc^2}\eta (x){T_{\lambda }(x)\wedge dd^c\Vert x\Vert^2
\over
\Vert x-z\Vert^2}=(\eta (x)\sum_{j=1}^2t_{\lambda }^{jj}(x))*K(z)\eqno(2.4.4)
$$
\noindent
where $K(z):={1\over \Vert z\Vert^2}$. $U_{\lambda }$ are bounded in $L^1(
\Delta^2)$ because $T_{\lambda }$ are so. From [Sk], p.376, we have that
$$
{\partial^2U_{\lambda }(z)\over \partial z_{\alpha }\partial\bar z_{\beta }}
= \eta (z)\cdot t_{\lambda }^{\alpha \bar\beta }(z) + \sum_{j=1}^2
{\partial K\over \partial z_{\alpha }}*({\partial \eta \over \partial\bar x
_{\beta }}t_{\lambda }^{j \bar j}-{\partial \eta \over \partial\bar x_j}
t_{\lambda }^{j \bar\beta }) +
$$

$$
 + \sum_{j=1}^2{\partial K\over \partial\bar z_j}*({\partial \eta \over
\partial x_{\alpha }}t_{\lambda }^{j \bar\beta }-{\partial \eta \over \partial
 x_j}t_{\lambda }^{\alpha \bar\beta }).\eqno(2.4.5)
$$
\noindent Using the fact that ${\sl supp}\nabla\eta \subset A^2(3/4,7/8)$ and
usual properties of convolutions, we see that the family $\{ dd^cU_{\lambda }\}
 $ is uniformly $C^{\infty }$-bounded on $A^2(1/2,3/4)$. So $\{ U_{\lambda }\}
$ are $C^{\infty }$-bounded on $A^2(1/2,3/4)$. Denote by $t_{\lambda \eps }^
{\alpha \bar\beta }$ the smoothing of $t_{\lambda }^{\alpha \bar\beta }$ by
convolution and by $S_{\lambda }^{\delta }$ the $\delta $-neighborhood of
$S_{\lambda }$. We have that
$$
\int_{\Delta^2(1/2)\setminus S_{\lambda }^{\delta }}(f^*w)^2 = \lim_{\eps
\searrow 0}\int_{\Delta^2(1/2)\setminus S_{\lambda }^{\delta }}T_{\lambda
\eps
}\wedge T_{\lambda \eps }\le \lim_{\eps \searrow 0}\int_{\Delta^2(1/2)}T_
{\lambda \eps }\wedge T_{\lambda \eps } =
$$
$$
 = \lim_{\eps \searrow 0}\int_{\Delta^2(1/2)}\eps_{\alpha \bar\beta \delta
\bar\gamma }t_
{\lambda \eps }^{\alpha \bar\gamma }\cdot t_{\lambda \eps }^{\delta \bar\beta
}d^4V =  \lim_{\eps \searrow 0}\int_{\Delta^2(1/2)}\eps_{\alpha \bar\beta
\delta \bar\gamma }
[{\partial^2U_{\lambda \eps }\over \partial z_{\alpha }\partial \bar z_{\gamma
 }}  -
$$
$$
 - \sum_{j=1}^2{\partial K_{\eps }\over \partial z
_{\alpha }}*({\partial \eta \over \partial\bar x_{\gamma }}t_{\lambda }^{j\bar
j} - {\partial \eta \over \partial\bar x_j}t_{\lambda }^{j \bar\gamma }) -
\sum_{j=1}^2{\partial K\over \partial\bar z_j}*({\partial \eta \over
\partial x_{\alpha }}t_{\lambda }^{j \bar\gamma }-{\partial \eta \over \partial
 x_j}t_{\lambda }^{\alpha \bar\gamma })]\cdot
$$
$$
 \cdot [{\partial^2U_{\lambda \eps }\over \partial z_{\delta }\partial \bar z_
{\beta }}  - \sum_{j=1}^2{\partial K_{\eps }\over \partial z
_{\delta }}*({\partial \eta \over \partial\bar x_{\beta }}t_{\lambda }^{j\bar
j} - {\partial \eta \over \partial\bar x_j}t_{\lambda }^{j \bar\beta }) -
\sum_{j=1}^2{\partial K\over \partial\bar z_j}*({\partial \eta \over
\partial x_{\delta }}t_{\lambda }^{j \bar\beta }-{\partial \eta \over \partial
 x_j}t_{\lambda }^{\delta \bar\beta })]d^4V =
$$
$$
 = \lim_{\eps \searrow }\int_{\partial \Delta^2(1/2)}J_{\lambda \eps } <
\infty ,
$$
\noindent
where by $J_{\lambda \eps }$ we denote the appropriate expression, which
remains after integrating. This $J_{\alpha \eps }$ is clearly $C^{
\infty }$-bounded on $A^2(1/2,3/4)$ uniformly on $\lambda $ and $\eps $.
\smallskip
\hfill{q.e.d.}
\smallskip
One has the following obvious
\smallskip\noindent\bf
Corollary 2.4.1. \it Let ${\cal F}$ be a family of meromorphic mappings
from $\Delta^2$ to a disk-convex K\"ahler space $X$, which is equicontinuous
on the Hartogs domain $H^2(r)\subset \Delta^2$. Then for any $\rho <1$ there
is a constant $C=C_{\rho , {\cal F}}$  such that ${\sl vol}(\Gamma_
{f_{\lambda }})\le C$ on $\Delta^2_{\rho }$, for all $f_{\lambda }\in
{\cal F}$.

\smallskip\rm
Remark that this statement doesn't follows from the Oka-type estimates for
the volumes of analytic sets of masses of currents, compare [FS-2]. We shall
need this estimate in the proof of Theorem 3 below.

In this regard we whant also to propose the following

\smallskip\noindent\sl
Conjecture. \it Let ${\cal F}$ be a family of meromorphic mappings
from $\Delta^n$ to a disk-convex K\"ahler space $X$, which is equicontinuous
on the Hartogs domain $H^n(r)\subset \Delta^n$. Then for any $\rho <1$ there
is a constant $C=C_{\rho , {\cal F}}$  such that ${\sl vol}(\Gamma_{
f_{\lambda }})\le C$ on $\Delta^n_{\rho }$, for all $f_{\lambda }\in {\cal F}$.
\smallskip\rm
Equicontinuity condition here, as well as in Corollary 2.4.1, means in
particular, that $f_{\lambda }$ are \it holomorphic \rm on $H^n(r)$.
Example of Shiffman and Taylor, see [Si-2], shows that one should essentially
use the fact that currents $T_{\lambda }$ are preimages of the K\"ahler form
by the mappings $f_{\lambda }$!

\smallskip\noindent\sl
2.5. Fatou sets of meromorphic self-maps of compact complex surfaces.
\smallskip\rm

We shall consider here a special case, when our family ${\cal F}$
is the family $ \{ f^n\} $ of iterates of some meromorphic self-map $f:X\to
X$ of an algebraic surface $X$, then strong Fatou set (i.e. set of convergency
of the family of iterates) coincides with weak Fatou set unless a quite
special case occurs. Namely we have the following

\bigskip\noindent\bf
Theorem 3. \it  Let $f$ be a meromorphic self-map of a compact K\"ahler
surface $X$. Denote by $\Phi_s $ the (strong) Fatou set of $f$ and by $\Phi_w$
 the weak Fatou set of $f$. Then:
\smallskip
\noindent
(i) $\Phi_w$ is a Levi-pseudoconvex open subset of $X$;
\smallskip
\noindent
(ii) If $\Phi_s $ is not equal to $\Phi_w$, then:


(a) $\Phi_w\supset X\setminus C$, where $C$ is a  rational curve in $X $;


(b) any weakly converging subsequence $\{ f^{n_k}\} $ converge strongly on

$X \setminus (C\cup \{  $ finite set $\} ) $,
and its  weak limit $f_\infty $  is a degenerate mapping of $X$ onto $C$.

\smallskip\rm

First let us start with a simple {\sl Lemma}. Let a meromorphic map
$f:\Omega \to X$ of complex spaces is given. By an image of a point $a$
(or more generally a set $A$) one understands $f[A]:=\{ x\in X: \exists
a\in A $ s.t. $(a,x)\in \Gamma_f\} $. If $A$ is an analytic subset, then
by $f\mid_A(A)$ we understand the image of $A$ under the restriction of
$f$ onto $A$. An analytic set $D(f):=p_1(\{ (z,x)\in \Omega \times X:
\dim_{(z,x)}p_2^{-1}(x)\ge 1\} )$ is called a locus of degeneration of $f$.
$f:\Omega \to X$ is degenerate if $D(f)=\Omega $.
\smallskip\noindent\bf
Lemma 2.5.1. \it Let $z\in \Phi_s $ (corr. in $\Phi_w$) and take some $l\ge 1$.
Then $f^l[z]\setminus f^l\mid_{D(f^l)}(D(f^l))\subset \Phi-s $ (corr. $\Phi_w$
).
\smallskip\noindent\sl
Proof. \rm Let  $v\in f^l[z]\setminus f^l\mid_{D(f^l)}(D(f^l))$, then there
are
 neighborhoods $U\ni z$ and $V\ni v$ such that  $f^{-l}:V\to U$ is a
multivalued \it holomorphic map. \rm Take some
sequence $\{ f^{n_k}\} \subset \{ f^n\} $. From sequence $\{ f^{n_k+l}\} $ by
assumption one can subtract a converging (corr. weakly converging) subsequence
 $\{ f^{n_{k_i}+l} \}$ on $U$. \rm So $f^{n_{k_i}}
 = f^{n_{k_i}+l}\circ f^{-l}$ will converge (corr. weakly converge) on
$V$.
\smallskip
\hfill{q.e.d.}
\smallskip

\smallskip\noindent\sl
Proof of Theorem 4. \rm (i) This is consequence of Corollary 1.2.1 and the
Hartogs-type extension theorem for maromorphic mappings into K\"ahler
manifolds, see [Iv-2].

\smallskip\noindent (ii) If $\Phi_s \not= \Phi_w$ then there exists a
point $p\in \Phi_w$, a ball $\bb $ centered at $p$, a subsequence of iterates
$\{ f^{n_k}\} $, which s-converge on $\bar \bb \setminus \{ p\} $ to a
meromorphic map $f_{\infty }:\bar\bb \to X$, but not converge on any
neighborhood of $p$. In particular this means that $p\in I(f_{\infty })$ by
Rouche principle.
Taking $\bb $ small enough, we can suppose that $p$ is the only fundamental
point of $f_{\infty }$ in $\bar\bb $.

By {\sl Theorem 2 } ${\sl vol}(\Gamma_{f^{n_k}})$ are uniformly bounded on
$\bar\bb
$. So $\{ \Gamma_{f^{n_k}}$ are converging (after going to subsequence) in
Hausdorff metric
to $\Gamma_{f_{\infty }}\cup (\{ p\} )\times X$. Put $C=f_{\infty }[p]$. $C$
is a finite union $\bigcup_{i=1}^NC_i$ of rational curves. Take a point $q\in
X\setminus C$. Then for $k\gg 1$ $q\in f^{n_k}(\bb \setminus \{ p\} )$. If
moreover $q\not\in f^{n_k}(D(f^{n_k}))$ then $q\in \Phi_w$ by Lemma 2.5.1.
But $\bigcup_k
f^{n_k}(D(f^{n_k}))$ is at most countable set of points and $\Phi_w$ is
Levi-pseudoconvex.
So $\Phi_w\supset X\setminus C$.

Take now a point $x\in C$. Suppose that $\Gamma_{f_{\infty }}\cap (X\times
\{ x\} )$ has $(p,x)$ as isolated point. Then we can find a neighborhoods
$W\ni p$ and $V\ni x$ with $(\partial W\times \bar V)\cap \Gamma_{f_{\infty }}
=\emptyset $. Thus $\partial W\times \bar V)\cap \Gamma_{f^{n_k}}=\emptyset $
for $k$ big enough. So by  $f^{n_k}(W)\supset V$. Thus $\Phi_w\supset C$.

Let us distinguish two cases.
\smallskip\noindent\sl
Case 1. \rm  $C$ has such a point $x$.
\smallskip By pseudoconvexity $\Phi_w$ contains a component of $C$, to which
$x$ belongs. By connectivity of $C$ $\Phi_w\supset C$ and thus  $\Phi_w=X$.
 In this case our sequence $\{ f^{n_k}\} $ strongly converges on $X\setminus
\{ p_1,...,p_N\} $.
From {\sl Theorem 2} we see that ${\sl vol}(\Gamma_{f^{n_k}})$ are uniformly
bounded but not less then $N\cdot {\sl vol}X$. This is possible only if
$f$ has degree one, $N=1$ and $f_{\infty }$ is degenerate map onto $C$. $C$
in this case should consist only from one component.
\smallskip\noindent\sl
Case 2. \rm For all points $x\in C$ ${\sl dim}_{(p,x)}\Gamma_{f_{\infty }}
\cap X\times \{ x\} >0$.
\smallskip
 Then $f_{\infty }$ is a degenerate mapping of $X\setminus C$ onto $C$. By
 the same
 reasoning as in the previous case $p$ is a single point in $\Phi_w
\setminus \Phi_s $.
\smallskip
\hfill{q.e.d}
\smallskip
The following example shows that the situation described in part (ii) of this
theorem can really happen. Let $X=\cc\pp^2$ and $f:[z_0:z_1:z_2]\to [z_0:
2z_1:2z_2]$. Then for this $f$ we have the phenomena described above with
$p=[1:0:0]$ and  $C=\{ z_0=0\} $.

\smallskip\noindent\bf
Remark. 1. \rm The statement (i) of this Theorem is valid for meromorphic
self-maps of compact K\"ahler manifolds of any dimension for the same reason.
It is also valid for all compact complex surfaces. This follows from
Corollary 2.3.3 and the fact that every compact K\"ahler surface carries
a pluriclosed metric form.

\noindent\bf
2. \rm In the firsrt Case, see the  proof, $f$ is a bimeromorhpic
automorphism of $X$. Probably one should expect that if $\Phi_s\not= \Phi_w$
then $f$ is nessessaliry  a bimeromorphic automorphism. The dynamics
of birational automorphisms of $\pp^2$ where recently studied in [D].

\magnification=\magstep1
\spaceskip=4pt plus3.5pt minus 1.5pt
\spaceskip=5pt plus4pt minus 2pt
\font\csc=cmcsc10
\font\tenmsb=msbm10

\def\cc{\hbox{\tenmsb C}}
\newdimen\length
\newdimen\lleftskip
\lleftskip=2.5\parindent
\length=\hsize \advance\length-\lleftskip
\def\entry#1#2#3#4\par{\parshape=2  0pt  \hsize%
\lleftskip \length%
\noindent\hbox to \lleftskip%
{\bf[#1]\hfill}{\csc{#2 }}{\sl{#3}}#4%
\medskip
}
\ifx \twelvebf\undefined \font\twelvebf=cmbx12\fi

\bigskip\bigskip
\bigskip\bigskip
\centerline{\twelvebf References.}
\bigskip


\entry{B-S}{Bedford E., Smillie J.:}{ Polynomial diffeomorphisms fo $\cc^2$
.}  Invent.\ Math. (1990) {\bf87}, 69-99.

\entry{C-T}{Cartan B., Thullen P.:}{Zur der Singularit\"aten der Funktionen
mehrerer Verandenlichen: Regularit\"ats und Konvergenzbereiche.} Math.\ Ann.
 {\bf106}, 617-647 (1932).


\entry{D}{Diller J.:}{Dynamics of birational maps of $\pp^2$.}
\ Preprint.

\entry{FS-1}{Fornaess J.-E., Sibony N.:}{Complex dynamics in higher
dimension II.} Modern Methods in Complex Analysis, Princeton Univ. Press
(1995), pp.135-182.

\entry{FS-2}{Fornaess J.-E., Sibony N.:}{Oka's unequality for currents and
applicatins.} to \ appear \ in \ Math.\ Ann .

\entry{Fj}{Fujimoto H.:}{On families of meromorphic maps into the complex
projective space.} Nagoya\ Math.\ J. {\bf54}, 21-51. (1974).

\entry{Ga}{Gauduchon P.:}{Les metriques standard d'une surface a premier
nombre de Betti  pair.} Asterisque.\ Soc.\ Math.\ France.  {\bf126}, 129-135,
(1985).






\entry{Iv-1}{Ivashkovich S.:}{Rational curves and Extensions of Holomorphic
mappings.} Proc.\ Symp.\ Pure \ Math. {\bf52} Part 1, 93-104, (1991).

\entry{Iv-2}{Ivashkovich S.:}{Spherical shells as obstructions for the
extension of  holomorphic mappings.} The\ Journal\ of\ Geometric\ Analysis.
{\bf2}, N 4, 351-371, (1992).

\entry{Iv-3}{Ivashkovich S.:}{The Hartogs-type extension theorem for the
meromorphic maps into compact K\"ahler manifolds.} Invent.\ math. {\bf109}
 , 47-54, (1992).

\entry{Iv-4}{Ivashkovich S.:}{Continuity principle and extension properties
of  meromorphic mappings with values in non K\"ahler manifolds.}
MSRI-Preprint, 1997-033.

\entry{Iv-5}{Ivashkovich S.:}{One example in concern with extension and
separate analyticity properties of meromorphic  mappings.} Preprint SFB -237,
 Bochum (1994), to appear in Amer. J. Math.

\entry{Ka}{Kato, M.:}{Compact complex manifolds containing "global
spherical shells".} Proc.\ Int.\ Symp.\ Alg.\ Geom.\ Kyoto. 45-84, (1977).






\entry{Re-1}{Remmert R.:}{Funktionentheorie 1} Springer (1995).

\entry{Re-2}{Remmert R.:}{Holomorphe und meromorphe Abbildungen komplexer
R\"aume.} Math.\ Ann. {\bf133}, 328-370, (1957).




\entry{Si-1}{Siu Y.-T.:}{Analyticity of Sets associated to Lelong Numbers and
the extension of Closed Positive Currents.} Invent.\ math. {\bf27}, 53-156,
(1974).

\entry{Si-2}{Siu Y.-T.:}{Extension of meromorphic maps into K\"ahler manifolds
.} Ann. \ Math. {\bf102}, 421-462, (1975).

\entry{Si-3}{Siu Y.-T.:}{Every stein subvariety admits a stein neighborhood.}
 Invent.\ Math. {\bf38}, N 1, 89-100, (1976).

\entry{Sk}{Skoda H.:}{Prolongement des courans, positifs, ferm\'es de masse
finie.} Invent.\ math. {\bf66}, 361-376, (1982).

\entry{St}{Stolzenberg G.:}{Volumes, Limits and Extension of Analytic
Varieties.} Springer \ Verlag (1966).

\bigskip\bigskip
Universite de LIlle-I

U.F.R. de Mathematiques

59655 Lille

France

ivachkov@gat.univ-lille1.fr

\bigskip\bigskip
IAPMM Acad Sci. of Ukraine

Naukova 3/b, 290053 Lviv

Ukraine
\end